\documentclass[a4paper,12pt]{elsarticle}
\journal{arXiv}
\usepackage{amsmath}
\usepackage[english]{babel}
\usepackage{graphicx}
\usepackage[usenames,dvipsnames]{color} 
\usepackage[numbers]{natbib}
\usepackage[footnotesize,bf,sf,center]{caption}
\usepackage{float}
\usepackage[dvips]{epsfig}
\usepackage[bottom=2.5cm,top=2.5cm,left=3cm,right=2cm]{geometry}
\usepackage{calc}
\usepackage{soul}
\usepackage{amsmath,amsfonts,amssymb}
\usepackage{subfigure}
\usepackage{setspace}
\usepackage{amssymb}
\usepackage{mathrsfs}
\usepackage{algorithm}
\usepackage{algpseudocode}
\usepackage{dsfont}
\usepackage{epstopdf}
\usepackage{subfig}%
\usepackage[utf8]{inputenc}
\usepackage{nicefrac,xfrac}
\newcommand{\ve}[1]{\mbox{\boldmath $#1$}}
\newcommand{\vesmall}[1]{\mbox{\boldmath \scriptsize $#1$}}
\newtheorem{theorem}{Theorem}[section]

\newtheorem{remark}{Remark}[section]

\newdefinition{rmk}{Remark}
\newcommand{\proof} [1]{ \noindent {\bf Proof.} #1 \hfill\rule{0.5em}{1.2ex} \par\medskip}
\usepackage{graphicx}
\begin{document}
	\begin{frontmatter}
		
		\renewcommand\arraystretch{1.0}
		
		\title{\textbf{Consistent splitting SAV schemes for finite element approximations of incompressible flows}}
		\author{
			{\bf Douglas R.~Q.~Pacheco} $^{1,2,3,}$
			\footnote[0]{{\sf Email address:} {\tt pacheco@ssd.rwth-aachen.de}, 
				{\sf corresponding author}}\\
			{\small ${}^{1}$ Chair for Computational Analysis of Technical Systems, RWTH Aachen University, Germany}\\
			{\small ${}^{2}$ Chair of Methods for Model-based Development in Computational Engineering, RWTH}\\
			{\small ${}^{3}$ Center for Simulation and Data Science (JARA-CSD), RWTH Aachen University, Germany}}
		\begin{keyword}
			Fractional-step methods \sep SAV schemes
			\sep Incompressible Navier--Stokes \sep IMEX methods  \sep Splitting schemes \sep BDF2
		\end{keyword}
		\begin{abstract}
			Consistent splitting schemes are among the most accurate pressure segregation methods, incurring no splitting errors or spurious boundary conditions. 
            Nevertheless, their theoretical properties are not yet fully understood, especially when finite elements are used for the spatial discretisation. This work proposes a simple scalar auxiliary variable (SAV) technique that, when combined with standard finite elements in space, guarantees unconditional stability for first- and second-order consistent splitting schemes. The framework is implicit-explicit (IMEX) and only requires solving linear transport equations and a pressure Poisson problem per time step. Furthermore, pressure stability is attained with respect to a stronger norm than in classical projection schemes, which allows eliminating the inf-sup compatibility requirement on the velocity-pressure pairs. The accuracy of the new framework is assessed through numerical examples.
		\end{abstract}
	\end{frontmatter}
	
	\section{Introduction}
Fractional-step methods, which allow computing velocity and pressure separately, make up for the majority of high-performance incompressible flow solvers. Since reviewing the history of such methods has become absolute commonplace in the literature, the reader is here simply referred to a overview articles \cite{Guermond2006,Badia2008}. What is more important here is to provide a brief \textsl{taxonomy} of fractional-step methods. The most common approaches can be grouped into three families: pressure-correction, velocity-correction and consistent splitting schemes. The first two of those families are derived from the spatially continuous level (although fully algebraic versions also exist \cite{Badia2008}) and are often called projection methods, since an intermediate velocity is projected onto a divergence-free space. Pressure-correction methods, the most popular among them, can normally not surpass or even reach full second-order convergence, due to inherent inconsistencies in their formulation \cite{Guermond2005}. Velocity-correction methods, on the other hand, can often achieve higher orders, but usually at the cost of more substeps \cite{Karniadakis1991}. The third family of segregation methods is based on an equivalent reformulation of the Navier--Stokes system at the continuous level, replacing the incompressibility constraint with a consistent pressure Poisson equation (PPE). This is the framework the present work addresses. 

Although a class of consistent splitting methods was already introduced in Ref.~\cite{Guermond2003}, the one addressed here was proposed shortly after \cite{Johnston2004,Liu2009}. In the words of \citet{Li2023SAV}, ``the consistent splitting scheme has outstanding advantages''. Indeed, it is free from splitting errors and hence arbitrarily high-order accurate (the full order of the temporal discretisation is attained). Furthermore, differently from, e.g., rotational pressure-correction methods, consistent splitting schemes do not require inf-sup-stable finite element pairs. In fact, while there is vast numerical evidence supporting that \cite{Guermond2006,Guermond2003,Wu2022}, the analysis is rather scarce. There are articles proving convergence in the context of spectral \cite{Johnston2004}, Legendre--Galerkin \cite{Wu2022}, marker-and-cell \cite{Li2023} and $C^1$ finite element methods \cite{Liu2009C1}, but so far nothing on standard $C^0$ elements \cite{Li2023}. In this context, the present article is, to the best of the author's knowledge, the first one to prove temporal stability for a consistent splitting scheme discretised in space with Lagrangian finite elements. To accomplish that, the scheme is combined with a simple SAV technique.

The \textsl{scalar auxiliary variable} method was first introduced for gradient flows \cite{Shen2019} but has since been used in a variety of fluid flow problems, especially in combination with fractional stepping \cite{Li2021}. Although SAV schemes have been typically employed to allow treating convection explicitly, some approaches have very recently been proposed focusing on stabilising pressure extrapolations \cite{Sun2024,Obbadi2025}. In the present work, the auxiliary variable is introduced to help control the explicit pressure term in the momentum equation. The resulting scheme uses an extrapolated pressure to update the velocity and the SAV, which are then fed into the consistent PPE for the pressure update. It will be shown that the second-order version of the scheme, when combined with any standard, $C^0$--continuous finite element pair, is temporally stable for any time-step size. The analysis is carried out for no-slip boundary conditions, but the method is presented in a more general boundary condition setting. The scheme is linearised and IMEX, in the sense that both the pressure and the convective velocity in the momentum equation are treated explicitly. 

 This article continues with the following structure. Sec.~\ref{sec_pre} introduces the PDE system and a consistent reformulation derived at the continuous level. The time-stepping scheme is proposed in Sec.~\ref{sec_schemes}, with the final algorithm presented in Subsec.~\ref{subsec_alg}. The stability analysis is carried out in Sec.~\ref{sec_analysis}, where unconditional temporal stability is proved for its second-order version. Sec.~\ref{sec_comparison} compares the new method to similar approaches with respect to stability properties. Sec.~\ref{sec_examples} then presents various numerical examples assessing accuracy and stability, while concluding remarks are drawn in Sec.~\ref{sec_Conclusion}.
	
\section{Preliminaries}\label{sec_pre}
\subsection{Notation and useful identities and inequalities}
Let $\Omega\subset\mathbb{R}^{d}$, $d=2$ or $3$, be a bounded Lipschitz domain. The $L^2(\omega)$ product of two functions is denoted here by $\langle \cdot ,\cdot \rangle_{\omega}$, with the domain omitted when $\omega = \Omega$. The $L^2(\Omega)$ and $L^{\infty}(\Omega)$ norms will be written as $\| \cdot \|$ and $\| \cdot \|_{\infty}$, respectively. Both the data and the solution of the continuous problem are assumed as sufficiently regular. Approximate or discrete values of quantities at different time steps will be denoted with subscripts: $\ve{u}_{n}$, for instance, denotes the approximation of $\ve{u}$ at the $n$th time step $t_n = n\tau$, with the time-step size $\tau$ assumed as constant, for simplicity. The following identity will be useful:
\begin{align}
     2\langle3\ve{v}_{n+1}-4\ve{v}_{n}+\ve{v}_{n-1},\ve{v}_{n+1}\rangle &= \|\ve{v}_{n+1}\|^2 - \|\ve{v}_{n}\|^2 + \|\ve{v}_{n+2}^{\star}\|^2 - \|\ve{v}_{n+1}^{\star}\|^2 +\|\delta^2\ve{v}_{n+1}\|^2  ,
    \label{identityBDF2}
\end{align}
which holds for any scalar or vector-valued $\ve{v}_{n+1},\ve{v}_{n},\ve{v}_{n-1}$, where
\begin{align*}
\ve{v}_{n+1}^{\star} &:= 2\ve{v}_{n}-\ve{v}_{n-1}\, ,\\
\delta^2\ve{v}_{n+1}&:=\ve{v}_{n+1}-2\ve{v}_{n}+\ve{v}_{n-1} = \ve{v}_{n+1} - \ve{v}_{n+1}^{\star} \, .
\end{align*}
	Moreover, provided that $\ve{v}\cdot\ve{n}=0$ on $\partial\Omega$, there holds
	\begin{align}
		\left\langle \ve{w}\cdot\nabla\ve{v} + \frac{\nabla\cdot \ve{w}}{2}\ve{v},\ve{v}\right\rangle &= 0
		\label{skewSymU}
	\end{align}
	for any $(\ve{v},\ve{w})$ sufficiently smooth.

We will consider shape-regular, globally quasi-uniform triangulations $\mathcal{T}$ of $\Omega$, with  
\begin{align}
1 \leq \frac{h}{h_e} \leq Q, \quad Q \geq 1\, ,\label{quasiUniform}
\end{align}
where $h_e$ denotes the size of element $\Omega_e$, and $h$ is the maximum element size in $\mathcal{T}$. Based on that, let $X_h$ be a standard Lagrangian finite element space. Then, the following trace inequality holds for any $v_h\in X_h$ (Lemma 1.46 in Ref.~\cite{DiPietro2012}):
\begin{align}
\sqrt{h_e}\| v_h \|_{L^2(F_e)} \leq C_{\text{tr}}\| v_h \|_{L^2(\Omega_e)} \, \ \text{for any face $F_e$ of $\Omega_e$,}\label{traceIneq}
\end{align}
with $C_{\text{tr}}$ depending only on shape-regularity constants, the spatial dimension $d$ and the polynomial degree of $X_h$.

\subsection{The incompressible Navier--Stokes system}
Considering a finite time interval $(0,T]$, the incompressible flow of a Newtonian fluid in $\Omega$ can be described by 
	\begin{align}
\nabla\cdot\ve{u} &= 0  && \text{in} \ \ \Omega\times(0,T]\, , \label{incompressibility}\\
\partial_t\ve{u} + \ve{u}\cdot\nabla\ve{u} - \nu\Delta\ve{u}  + \nabla p &= \ve{f} && \text{in} \ \ \Omega\times(0,T]\, ,\label{momentum}
\end{align}
where the unknowns are the velocity $\ve{u}$ and the pressure $p$, while the viscosity $\nu>0$ and the force $\ve{f}$ are given. The remaining problem data are
\begin{align}
\ve{u} &= \ve{u}_0 && \text{at} \ \ t=0,\label{initialCondition}\\
\ve{u} &= \ve{g} && \text{on} \ \ \Gamma_D\times(0,T]\, ,\label{DirichletBC}\\
(\nu\nabla\ve{u})\ve{n} - p\ve{n} &= \ve{t} && \text{on} \ \ \Gamma_N\times(0,T]\, , \label{tractionBC}
\end{align}
in which $\ve{n}$ is the outward unit normal vector, ($\ve{u}_0,\ve{g},\ve{t}$) are known vectors, and ($\Gamma_D,\Gamma_N$) form a non-overlapping decomposition of the boundary $\Gamma:=\partial\Omega$.

\subsection{Equivalent reformulation}
Before the temporal discretisation is introduced, the Navier--Stokes system will be recast into an equivalent form. First, let us define the auxiliary variable $\psi(t) \equiv 1$, which will be treated as an additional unknown. Due to \eqref{incompressibility}, the function $\psi$ satisfies, for any $\alpha\in\mathbb{R}$, the initial value problem
\begin{align*}
    \psi'(t) = -\alpha\underbrace{\int_{\Omega} p\nabla\cdot\ve{u}\, \mathrm{d}\Omega}_{=\, 0}\, , \ \ \psi(0)=1\, .
\end{align*}
Based on that, the original equations \eqref{incompressibility}--\eqref{tractionBC} can be replaced by the equivalent system
\begin{align}
    \partial_t\ve{u} + \ve{u}\cdot\nabla\ve{u} + \frac{\nabla\cdot\ve{u}}{2}\ve{u} - \nu\Delta\ve{u} -\gamma\nu\nabla(\nabla\cdot\ve{u})  + \psi\nabla p &= \ve{f} && \text{in} \ \, \Omega\times(0,T]\, ,\label{momentumModified}\\
    \psi\Delta p  + \nabla\cdot\left[\nu\nabla\times(\nabla\times\ve{u}) + \ve{u}\cdot\nabla\ve{u} + \frac{\nabla\cdot\ve{u}}{2}\ve{u} - \ve{f}\right]\label{PPE} &= 0 && \text{in} \ \, \Omega\times[0,T]\, ,\\
       \frac{\mathrm{d}\psi}{\mathrm{d}t} +\alpha\big\langle p,\nabla\cdot\ve{u}\big\rangle &= 0 && \text{in} \ \, (0,T]\, ,\label{ODE}
\end{align}
with initial and boundary data
\begin{align}
    &\psi = 1 && \text{at} \ \, t=0 \label{psi0}\\
 &\ve{u} = \ve{u}_0 && \text{at} \ \, t=0\\
& \ve{u} = \ve{g} && \text{on} \ \ \Gamma_D\times(0,T]\, ,\\
&(\nu\nabla\ve{u})\ve{n} + (\gamma\nu\nabla\cdot\ve{u} - \psi p)\ve{n} = \ve{t} && \text{on} \ \ \Gamma_N\times(0,T]\, , \label{BCwithGradDiv} \\
&\psi p = \nu\nabla\ve{u}:(\ve{n}\otimes\ve{n})  - \ve{t}\cdot\ve{n}&& \text{on} \ \Gamma_N\times[0,T]\, , \label{PPE_DBC}\\
&\psi\partial_{\vesmall{n}} p = \ve{n}\cdot\left[\ve{f}-\partial_t\ve{g} - \ve{u}\cdot\nabla\ve{u} - \frac{\nabla\cdot\ve{u}}{2}\ve{u} -\nu\nabla\times(\nabla\times\ve{u}) \right] &&  \text{on} \ \Gamma_{D}\times[0,T]\, , \label{PPE_NBC}
\end{align}
where $\alpha$ and $\gamma$ are two positive, user-defined parameters. The motivation for introducing the SAV and the auxiliary equation \eqref{ODE} is to help control the pressure term in the momentum equation, which, in the case of consistent splitting schemes (differently from projection methods), does not have a balancing counterpart in the pressure equation.   

Showing that the solution $(\ve{u},p)$ of \eqref{incompressibility}--\eqref{tractionBC} also fulfils \eqref{momentumModified}--\eqref{PPE_NBC} is relatively simple. Using that $\nabla\cdot\ve{u}=0$, $\psi=1$, and 
\begin{align*}
    [(\nabla\ve{u})\ve{n}]\cdot\ve{n} &\equiv \nabla\ve{u}:(\ve{n}\otimes\ve{n})\, ,\\
    \Delta\ve{u} &\equiv \nabla(\nabla\cdot\ve{u}) - \nabla\times(\nabla\times\ve{u})\, , \\
   \nabla\cdot(\partial_t\ve{u}) &\equiv \partial_t(\nabla\cdot\ve{u})\, ,
\end{align*}
we can derive the the Poisson equation \eqref{PPE} and the Neumann condition \eqref{PPE_NBC} from \eqref{momentum} and \eqref{DirichletBC}, while the Dirichlet condition \eqref{PPE_DBC} comes directly from \eqref{tractionBC}. Proving the other direction, on the other hand, is not as straightforward, since in that case we are not allowed to assume incompressibility. Let us start by applying $\nabla\cdot$ to both sides of Eq.~\eqref{momentumModified} and adding the result to Eq.~\eqref{PPE}; most terms will cancel out, and what remains is
\begin{align}
    \partial_t(\nabla\cdot\ve{u}) -(1+\gamma)\nu\Delta(\nabla\cdot\ve{u}) = 0\, ,
    \label{heat}
\end{align}
that is, a heat equation for $\nabla\cdot\ve{u}$. To show that Eq.~\eqref{heat} admits only the trivial solution $\nabla\cdot\ve{u}\equiv 0$, we need boundary and initial conditions. Restricting Eq.~\eqref{momentumModified} to $\Gamma_D$, dotting the result with $\ve{n}$, and subtracting \eqref{PPE_NBC} gives us the Neumann condition
\begin{align}
    (1+\gamma)\nu\partial_{\vesmall{n}}(\nabla\cdot\ve{u}) = 0 \ \ \text{on} \ \, \Gamma_D\, .
\end{align}
Similarly, dotting \eqref{BCwithGradDiv} with $\ve{n}$ and subtracting the result from \eqref{PPE_DBC} yields the Dirichlet condition
\begin{align}
    \gamma\nu\nabla\cdot\ve{u} = 0 \ \ \text{on} \ \, \Gamma_N\, .\label{divU_DBC}
\end{align}
Therefore, with $\gamma>0$ and provided that $\nabla\cdot\ve{u}_0 = 0$, the unique solution of \eqref{heat}--\eqref{divU_DBC} is $\nabla\cdot\ve{u} \equiv 0$. Then, we have from \eqref{ODE}--\eqref{psi0} that $\psi(t) \equiv 1$. With that and incompressibility proved, we can reverse all the steps used to derive the modified system, which shows that the equivalence indeed holds. Some important remarks are:
\begin{itemize}
\item The equivalence showed above assumes that the solution $(\ve{u},p)$ is regular enough to allow differentiation (details on that can be found in the literature \cite{Heywood1982,Sani2006}). At the discrete level, we will use integration by parts to reduce the regularity requirements so that standard $H^1$--conforming elements can be used. 
    \item When discretised, the term $-\nu\gamma\nabla(\nabla\cdot\ve{u})$ will produce the so-called grad-div stabilisation \cite{John2016}. This technique improves mass conservation, which is important here because, at the discrete level, we will no longer have divergence-free velocities. 
    \item The term $(\nicefrac{1}{2}\nabla\cdot\ve{u})\ve{u}$ was added to the momentum equation to allow using identity \eqref{skewSymU}, which is essential to guarantee unconditional stability.
    \item As will be shown later in the analysis, the coefficient $\alpha>0$ in \eqref{ODE} could, in principle, be taken arbitrarily small. Setting $\alpha=\gamma=0$ would essentially recover the original consistent splitting by \citet{Liu2009}.
\end{itemize}

\section{Numerical method}\label{sec_schemes}
	\subsection{Weak PPE}
	Before introducing a discretisation, we must derive for the PPE a weak form free of second-order derivatives. Let us start by testing \eqref{PPE} with a function $q\in H^1(\Omega)$ such that $q=0$ on $\Gamma_N$. Then, integration by parts combined with \eqref{PPE_DBC} and \eqref{PPE_NBC} yields
\begin{align*}
        \psi\langle \nabla p,\nabla q\rangle = \left\langle \ve{f} - \ve{u}\cdot\nabla\ve{u} - \frac{\nabla\cdot\ve{u}}{2}\ve{u} - \nu\nabla\times(\nabla\times\ve{u})
 ,\nabla q \right\rangle - \left\langle\partial_t(\ve{n}\cdot\ve{g}),q\right\rangle_{\Gamma_D}.
 \end{align*}
The remaining second-order derivatives can be eliminated thanks to the identity \cite{Johnston2004}
\begin{align*}
-\left\langle \nabla\times(\nabla\times\ve{u}),\nabla q\right\rangle = \left\langle\nabla\times\ve{u}, \ve{n}\times\nabla q\right\rangle_{\Gamma}\, .
\end{align*}
That reduces the weak PPE to 
\begin{align}
        \psi\langle \nabla p,\nabla q\rangle = \left\langle \ve{f} - \ve{u}\cdot\nabla\ve{u} - \frac{\nabla\cdot\ve{u}}{2}\ve{u} 
 ,\nabla q \right\rangle + \left\langle\nu\nabla\times\ve{u}, \ve{n}\times\nabla q\right\rangle_{\Gamma_N} - \left\langle\partial_t(\ve{n}\cdot\ve{g}),q\right\rangle_{\Gamma_D} ,\label{weakPPE}
 \end{align}
where only first-order derivatives remain. 

\begin{remark}
    For variable $\nu$, the viscous part of the weak PPE \eqref{weakPPE} is generalised to
    \begin{align*}
       \left\langle\nu\nabla\times\ve{u}, \ve{n}\times\nabla q\right\rangle_{\Gamma_N} + \langle 2\nabla^{\top}\ve{u}\nabla\nu,\nabla q\rangle\, ,
    \end{align*}
see the derivation by \citet{Pacheco2022}. Also in that reference is a version of the PPE for non-constant density, as arising, for instance, in two-phase flows. 
\end{remark}

The time-stepping process in consistent splitting methods is simple: at each time step, we use an extrapolated pressure in the momentum equation to update the velocity, then feed that velocity into the PPE \eqref{weakPPE} to update the pressure. With the addition of the SAV, however, some key algorithmic modifications need to be introduced.

\subsection{Temporal discretisation}
Due to the full consistency of the present framework, arbitrarily high-order schemes can be constructed. Stability, on the other hand, depends on the specific time-stepping method used. We shall restrict the presentation to a second-order IMEX method whose stability is unconditional. Higher-order schemes are subject to the usual technical difficulties, such as the second Dahlquist barrier \cite{Dahlquist1963}. Nonetheless, stability may be provable for methods based on backward differentiation formulas (BDFs) of order 3--5 using specific techniques found in the literature \cite{Baker1982,Liu2013}.

To construct the fractional-step method, both $\partial_t\ve{u}$ and $\psi'(t)$ will be approximated by a second-order difference (BDF2). In the momentum equation, the convective velocity and the pressure will be extrapolated accordingly. Let us also consider two $H^1(\Omega)$--conforming finite element spaces $X_h$ and $Y_h$ for pressure and velocity, respectively. Then, each time step requires finding $(\ve{u}_{n+1},p_{n+1},\psi_{n+1}) \in X_h\times Y_h\times\mathbb{R}$ fulfilling $\ve{u}_{n+1}|_{\Gamma_D} = \ve{g}_{n+1}$ and $\psi_{n+1}p_{n+1}|_{\Gamma_N}=\mathcal{P}[\nu\nabla\ve{u}_{n+1}:(\ve{n}\otimes\ve{n})-\ve{t}_{n+1}\cdot\ve{n}]$ such that
\begin{align}
    &\left\langle \frac{3\ve{u}_{n+1} - 4\ve{u}_{n}+\ve{u}_{n-1}}{2\tau},\ve{v}\right\rangle + \left\langle  \ve{u}_{n+1}^{\star}\cdot\nabla\ve{u}_{n+1} + \frac{\nabla\cdot\ve{u}_{n+1}^{\star}}{2}\ve{u}_{n+1},\ve{v}\right\rangle + \nu\langle\nabla\ve{u}_{n+1},\nabla\ve{v} \rangle \nonumber  \\   &\quad +\, \gamma\nu\langle\nabla\cdot\ve{u}_{n+1},\nabla\cdot\ve{v}\rangle = \psi_{n+1}\langle p_{n+1}^{\star},\nabla\cdot\ve{v}\rangle + \left\langle\ve{f}_{n+1},\ve{v}\right\rangle + \langle\ve{t}_{n+1},\ve{v}\rangle_{\Gamma_N} ,\label{BDF2momentumWeak}\\
    &\frac{3\psi_{n+1}-4\psi_n+\psi_{n-1}}{2\tau} = -\alpha\langle p_{n+1}^{\star},\nabla\cdot\ve{u}_{n+1}\rangle\, ,\label{BDF2psi}\\
    &\psi_{n+1}\langle\nabla p_{n+1},\nabla q\rangle = \left\langle \ve{f}_{n+1} - \ve{u}_{n+1}\cdot\nabla\ve{u}_{n+1} - \frac{\nabla\cdot\ve{u}_{n+1}}{2}\ve{u}_{n+1}
 ,\nabla q \right\rangle\nonumber\\
 &\quad +\left\langle\nu\nabla\times\ve{u}_{n+1}, \ve{n}\times\nabla q\right\rangle_{\Gamma_N} - \left\langle\partial_t(\ve{n}\cdot\ve{g}_{n+1}),q\right\rangle_{\Gamma_D}\label{BDF2PPEweak}
\end{align}
for all $(\ve{v},q) \in X_h\times Y_h$ fulfilling $\ve{v}|_{\Gamma_D} = \ve{0}$ and $q|_{\Gamma_N}=0$, where
\begin{align*}
    \ve{u}^{\star}_{n+1} = 2\ve{u}_{n} - \ve{u}_{n-1} \ \ \text{and} \ \
    p^{\star}_{n+1} = 2p_{n} - p_{n-1}
\end{align*}
are formally second-order extrapolations, and $\mathcal{P}$ is a generic operator that makes the quantity $\nu\nabla\ve{u}_{n+1}(\ve{n}\otimes\ve{n})-\ve{t}_{n+1}\cdot\ve{n}$ continuous at the pressure nodes on $\Gamma_N$ (e.g., a simple averaging or projection operator). 

\subsection{Solution algorithm}\label{subsec_alg}
The term \textsl{scalar} auxiliary variable can be misleading: $\psi$ is indeed scalar, but not \textsl{just} in the same sense as the pressure is scalar; the extra unknown $\psi_{n+1}$ is, in fact, a \textsl{constant}, so that the discrete SAV system has only \textsl{one} additional degree of freedom in comparison to the original discrete problem. It does, however, introduce a nonlinearity in an otherwise sequential scheme. To circumvent that and construct a more efficient algorithm, let us take inspiration from an idea proposed by \citet{Li2021} and introduce the splitting 
\begin{align}
    \ve{u}_{n+1} = \mathbf{u}^1_{n+1} + \psi_{n+1}\mathbf{u}^2_{n+1}\, ,\label{splitting}
\end{align}
where $\ve{u}^1_{n+1}$ and $\ve{u}^2_{n+1}$ solve the (weak form of the) linear auxiliary problems
\begin{flalign*}
\begin{cases}
 \frac{1}{2\tau}(3\mathbf{u}_{n+1}^1-4\ve{u}_{n}+\ve{u}_{n-1}) + \ve{u}_{n+1}^{\star}\cdot\nabla\mathbf{u}_{n+1}^1 + \frac{1}{2}(\nabla\cdot\ve{u}_{n+1}^{\star})\mathbf{u}_{n+1}^1 - \nu\Delta\mathbf{u}_{n+1}^1  -\gamma\nu\nabla(\nabla\cdot\mathbf{u}_{n+1}^1) = \ve{f}_{n+1} \\
 \mathbf{u}_{n+1}^1|_{\Gamma_N} = \ve{g}_{n+1} \\
 [(\nu\nabla\mathbf{u}_{n+1}^1)\ve{n} + \gamma\nu\nabla\cdot\mathbf{u}_{n+1}^1]\big|_{\Gamma_D} = \ve{t}_{n+1}
\end{cases}&&
\end{flalign*}
and
\begin{flalign*}
\begin{cases}
 \frac{3}{2\tau}\mathbf{u}_{n+1}^2 + \ve{u}_{n+1}^{\star}\cdot\nabla\mathbf{u}_{n+1}^2 + \frac{1}{2}(\nabla\cdot\ve{u}_{n+1}^{\star})\mathbf{u}_{n+1}^2 - \nu\Delta\mathbf{u}_{n+1}^2  -\gamma\nu\nabla(\nabla\cdot\mathbf{u}_{n+1}^2) = -\nabla p_{n+1}^{\star} \\
 \mathbf{u}_{n+1}^2|_{\Gamma_D} = \ve{0} \\
 \big[(\nu\nabla\mathbf{u}_{n+1}^2)\ve{n} + \gamma\nu\nabla\cdot\mathbf{u}_{n+1}^2 - p^{\star}_{n+1}\ve{n}\big]\big|_{\Gamma_N} = \ve{0}
\end{cases}&&
\end{flalign*}
Notice that this splitting is \textsl{not} an approximation, since both the PDEs and the boundary conditions follow the linear combination \eqref{splitting}, thereby satisfying the superposition principle. After computing $\mathbf{u}_{n+1}^1$ and $\mathbf{u}_{n+1}^2$, we insert them into
\begin{align*}
    3\psi_{n+1}-4\psi_n+\psi_{n-1} &= - 2\tau\alpha\left\langle p_{n+1}^{\star},\nabla\cdot\mathbf{u}_{n+1}^1 + \psi_{n+1}\nabla\cdot\mathbf{u}_{n+1}^2\right\rangle\, , 
\end{align*}
which can be straightforwardly solved for $\psi_{n+1}$. The algorithm is then as follows:
	\begin{itemize}
		\item \textbf{Step 0:} To initialise, set $\psi_0=\psi_1=1$, compute $p_0$ from the initial data via \eqref{weakPPE}, and use a one-step method (e.g., BDF1 or Crank--Nicolson) to compute $\ve{u}_1$ and $p_1$. 
		
		\item \textbf{Step 1:} Find $\mathbf{u}_{n+1}^1\in X_h$, fulfilling $\mathbf{u}_{n+1}^1|_{\Gamma_D} = \ve{g}_{n+1}$, such that 
		\begin{align}
		   & \left\langle \frac{3\mathbf{u}_{n+1}^1 - 4\ve{u}_{n}+\ve{u}_{n-1}}{2\tau},\ve{v}\right\rangle + \left\langle  \ve{u}_{n+1}^{\star}\cdot\nabla\mathbf{u}_{n+1}^1 + \frac{\nabla\cdot\ve{u}_{n+1}^{\star}}{2}\mathbf{u}_{n+1}^1,\ve{v}\right\rangle\nonumber \\
            &+ \nu\langle\nabla\mathbf{u}_{n+1}^1,\nabla\ve{v} \rangle  +\, \gamma\nu\langle\nabla\cdot\mathbf{u}_{n+1}^1,\nabla\cdot\ve{v}\rangle = 
            \left\langle\ve{f}_{n+1},\ve{v}\right\rangle + \langle\ve{t}_{n+1},\ve{v}\rangle_{\Gamma_N}
		\end{align}
		for all $\ve{v}\in X_h$ fulfilling $\ve{v}|_{\Gamma_D} = \ve{0}$.
		\item \textbf{Step 2:} Find $\mathbf{u}_{n+1}^2\in X_h$, fulfilling $\mathbf{u}_{n+1}^2|_{\Gamma_D} = \ve{0}$, such that 
		\begin{align}
		   & \left\langle \frac{3\mathbf{u}_{n+1}^2}{2\tau},\mathbf{v}\right\rangle + \left\langle  \ve{u}_{n+1}^{\star}\cdot\nabla\mathbf{u}_{n+1}^2 + \frac{\nabla\cdot\ve{u}_{n+1}^{\star}}{2}\mathbf{u}_{n+1}^2,\mathbf{v}\right\rangle\nonumber \\
            &+ \nu\langle\nabla\mathbf{u}_{n+1}^2,\nabla\mathbf{v} \rangle  +\, \gamma\nu\langle\nabla\cdot\mathbf{u}_{n+1}^2,\nabla\cdot\mathbf{v}\rangle = 
            \left\langle p_{n+1}^{\star},\mathbf{v}\right\rangle 
		\end{align}
		for all $\mathbf{v}\in X_h$ fulfilling $\mathbf{v}|_{\Gamma_D} = \ve{0}$.
		
		\item \textbf{Step 3:} Update the SAV via
        \begin{align*}
            \psi_{n+1} = \frac{4\psi_{n} - \psi_{n-1} - 2\tau\alpha\langle p_{n+1}^{\star},\nabla\cdot\mathbf{u}^1_{n+1}\rangle}{3+2\tau\alpha\langle p_{n+1}^{\star},\nabla\cdot\mathbf{u}^2_{n+1}\rangle}\, ,
        \end{align*}
		then update the velocity $\ve{u}_{n+1} = \mathbf{u}^1_{n+1} + \psi_{n+1}\mathbf{u}^2_{n+1}$.

		\item \textbf{Step 4:} Find $p_{n+1}\in X_h$, fulfilling $\psi_{n+1}p_{n+1}|_{\Gamma_N}=\mathcal{P}[\nu\nabla\ve{u}_{n+1}:(\ve{n}\otimes\ve{n})-\ve{t}_{n+1}\cdot\ve{n}]$, such that
		\begin{align}
		    \psi_{n+1}\langle\nabla p_{n+1},\nabla q\rangle &= \left\langle \ve{f}_{n+1} - \ve{u}_{n+1}\cdot\nabla\ve{u}_{n+1} - \frac{\nabla\cdot\ve{u}_{n+1}}{2}\ve{u}_{n+1}
 ,\nabla q \right\rangle\nonumber\\
 &+\left\langle\nu\nabla\times\ve{u}_{n+1}, \ve{n}\times\nabla q\right\rangle_{\Gamma_N} - \left\langle\partial_t(\ve{n}\cdot\ve{g}_{n+1}),q\right\rangle_{\Gamma_D}\vphantom{\left\langle \ve{f}_{n+1} - \ve{u}_{n+1}\cdot\nabla\ve{u}_{n+1} - \frac{\nabla\cdot\ve{u}_{n+1}}{2}\ve{u}_{n+1}
 ,\nabla q \right\rangle}
		\end{align}
for all $q \in X_h$ fulfilling $q|_{\Gamma_N}=0$.
\end{itemize} 

The algorithm proposed above can be implemented very efficiently, as it only requires solving two linear advection-diffusion-reaction problems (both with the same coefficient matrix) and a Poisson equation. The temporal accuracy can also be increased by using higher-order BDFs for $\partial_t\ve{u}$ and $\psi'(t)$, combined with corresponding extrapolation formulas to define $\ve{u}_{n+1}^{\star}$ and $p_{n+1}^{\star}$.


\section{Temporal stability analysis}\label{sec_analysis}
Although the present fractional-step framework is derived for a general boundary condition setting, the analysis will (as usual) consider the homogeneous Dirichlet case. For concision, the case $\ve{f}=\ve{0}$ will be analysed, although the forcing term can be easily handled using standard techniques. 
	
\begin{theorem}[Temporal stability of the second-order scheme]
For any time-step size $\tau=T/N$, the fractional-step scheme \eqref{BDF2momentumWeak}--\eqref{BDF2PPEweak}, with $\ve{f}=\ve{0}$, satisfies 
\begin{align}
    \Phi_{N} +  \sum_{n=2}^N\left[\bigg(\frac{\delta^2\psi_{n}}{\sqrt{\alpha}}\bigg)^2 +  \|\delta^2\ve{u}_{n}\|^2 + 4\tau\nu\|\nabla\ve{u}_{n}\|^2 + 4\gamma\tau\nu\|\nabla\cdot\ve{u}_{n}\|^2\right]&= \Phi_{1} \, ,\label{stabilityU}\\
     \Phi_{N}  +
     \sum_{n=2}^{N}\left(\frac{8}{7}\tau\beta_{n}\psi_{n}^2\|\nabla p_{n}\|^2 + \bigg(\frac{\delta^2\psi_{n}}{\sqrt{\alpha}}\bigg)^2 +  \|\delta^2\ve{u}_{n}\|^2 + 4\tau\nu\gamma\|\nabla\cdot\ve{u}_{n}\|^2\right) &\leq  \Phi_{1} \, ,\label{stabilityP}
\end{align}
\end{theorem}
where
\begin{align}
    \Phi_{n} &:= \bigg(\frac{\psi_{n}}{\sqrt{\alpha}}\bigg)^2 + \bigg(\frac{\psi_{n+1}^{\star}}{\sqrt{\alpha}}\bigg)^2 + \|\ve{u}_{n}\|^2 + \|\ve{u}^{\star}_{n+1}\|^2\, ,\label{Phi}\\
    \beta_{n} &:= \left(\frac{1}{\nu}\|\ve{u}_{n}\|_{\infty}^2 + \frac{C\nu}{ h^2}\right)^{-1}, \label{beta}
\end{align}
and $C=(C_{\mathrm{tr}}^2Q)^2$ according to \eqref{quasiUniform} and \eqref{traceIneq}.

\proof{We start by setting $\ve{v}=4\tau\ve{u}_{n+1}$ in \eqref{BDF2momentumWeak} and using \eqref{identityBDF2}, which yields
\begin{align*}
    \|\ve{u}_{n+1}\|^2 - \|\ve{u}_{n}\|^2 + \|\ve{u}^{\star}_{n+2}\|^2 -\|\ve{u}^{\star}_{n+1}\|^2  +\|\delta^2\ve{u}_{n+1}\|^2 + 4\tau\nu\|\nabla\ve{u}_{n+1}\|^2 + 4\gamma\tau\nu\|\nabla\cdot\ve{u}_{n+1}\|^2\\
    = 4\tau\psi_{n+1}\big\langle p^{\star}_{n+1},\nabla\cdot\ve{u}_{n+1}\big\rangle ,
\end{align*}
where the convective term vanished due to identity \eqref{skewSymU}. Similarly, we multiply \eqref{BDF2psi} by $\alpha^{-1}\psi_{n+1}$ and add the result to the equation above, so that
\begin{align}
    &\Phi_{n+1} - \Phi_{n}  + \alpha^{-1}(\delta^2\psi_{n+1})^2 +  \|\delta^2\ve{u}_{n+1}\|^2 + 4\tau\nu\|\nabla\ve{u}_{n+1}\|^2 + 4\gamma\tau\nu\|\nabla\cdot\ve{u}_{n+1}\|^2= 0 \, ,\label{est1momentum}
\end{align}
with $\Phi_{n}$ defined in \eqref{Phi}. Adding Eq.~\eqref{est1momentum} from $n=1$ to $n=N-1$ proves \eqref{stabilityU}. Notice that Eq.~\eqref{est1momentum} already guarantees stability for the velocity, independently from the pressure, which is a remarkable property made possible by the SAV technique. 

To prove pressure stability, we set $q=2\psi_{n+1} p_{n+1}$ in \eqref{BDF2PPEweak}, so that
\begin{align}
    &2\psi_{n+1}^2\|\nabla p_{n+1}\|^2 = \nonumber\\
    &2\psi_{n+1}\left\langle\nu\nabla\times\ve{u}_{n+1}, \ve{n}\times\nabla p_{n+1}\right\rangle_{\Gamma} - 2\psi_{n+1}\left\langle\ve{u}_{n+1}\cdot\nabla\ve{u}_{n+1} + \frac{\nabla\cdot\ve{u}_{n+1}}{2}\ve{u}_{n+1}
 ,\nabla p_{n+1} \right\rangle . \label{est1p}
\end{align}
Due to the Hölder and Young inequalities, the convective term yields
\begin{align}
    &- \psi_{n+1}\left\langle2\ve{u}_{n+1}\cdot\nabla\ve{u}_{n+1} + (\nabla\cdot\ve{u}_{n+1})\ve{u}_{n+1}
 ,\nabla p_{n+1} \right\rangle\nonumber\\
 &\leq |\psi_{n+1}|\|\nabla p_{n+1}\|\|\ve{u}_{n+1}\|_{\infty}\left(2\|\nabla\ve{u}_{n+1}\|+\|\nabla\cdot\ve{u}_{n+1}\|\right)\nonumber\\
 &\leq \frac{4}{7}\psi_{n+1}^2\|\nabla p_{n+1}\|^2 + \frac{7}{4}\|\ve{u}_{n+1}\|_{\infty}^2\|\nabla\ve{u}_{n+1}\|^2 + \frac{1}{7}\psi_{n+1}^2\|\nabla p_{n+1}\|^2 + \frac{7}{4}\|\ve{u}_{n+1}\|_{\infty}^2\|\nabla\cdot\ve{u}_{n+1}\|^2\nonumber\\
& \leq \frac{5}{7}\psi_{n+1}^2\|\nabla p_{n+1}\|^2 + \frac{7}{2}\|\ve{u}_{n+1}\|_{\infty}^2\|\nabla\ve{u}_{n+1}\|^2 \ ,\label{estConvective}
\end{align}
where we have used that $\|\nabla\cdot\ve{v}\| \leq \|\nabla\ve{v}\|$ for any $\ve{v}\in \big[H^1_0(\Omega)\big]^d$. To estimate the boundary term, we can use the trace inequality \eqref{traceIneq} to write
\begin{align*}
\langle \nu\nabla\times\ve{u}_{n+1}, \ve{n}\times\nabla p_{n+1}\rangle_{\Gamma} &= \nu\sum_{\Omega_e\subset\mathcal{T}}\langle\nabla p_{n+1}\times\ve{n}, \nu\nabla\times\ve{u}_{n+1}\rangle_{\Gamma\cap\partial\Omega_e}\\
 &\leq  \nu\sum_{\Omega_e\subset\mathcal{T}}\| \nabla p_{n+1} \times\ve{n}\|_{L^2(\Gamma\cap\partial\Omega_e)}\|\nabla\times\ve{u}_{n+1} \|_{L^2(\Gamma\cap\partial\Omega_e)} \\
&\leq \nu\sum_{\Omega_e\subset\mathcal{T}}\| \nabla p_{n+1} \|_{L^2(\Gamma\cap\partial\Omega_e)} \| \nabla \ve{u}_{n+1} \|_{L^2(\Gamma\cap\partial\Omega_e)} \\
 &\leq C_{\text{tr}}^2\nu\sum_{\Omega_e\subset\mathcal{T}}{h_e}^{-1}\| \nabla p_{n+1} \|_{L^2(\Omega_e)} \| \nabla \ve{u}_{n+1} \|_{L^2(\Omega_e)} \, .
\end{align*}
That, combined with \eqref{quasiUniform} and the triangle and Young inequalities, yields
\begin{align}
2\psi_{n+1}\langle \nu\nabla\times\ve{u}_{n+1}, \ve{n}\times\nabla p_{n+1}\rangle_{\Gamma} &\leq   2QC_{\text{tr}}^2 h^{-1}\nu\| \nabla p_{n+1} \| \, \| \nabla \ve{u}_{n+1} \| \nonumber\\
&\leq \frac{2}{7}\psi_{n+1}^2\|\nabla p_{n+1} \|^2 + \frac{7C\nu}{ 2h^2}\nu\|\nabla\ve{u}_{n+1}\|^2 \, ,\label{estBoundary}
\end{align}
where $C= (C_{\text{tr}}^2Q)^2$. Inserting \eqref{estBoundary} and \eqref{estConvective} into \eqref{est1p} gives
\begin{align*}
    \psi_{n+1}^2\|\nabla p_{n+1}\|^2 
     \leq  \frac{7}{2}\left(\frac{1}{\nu}\|\ve{u}_{n+1}\|_{\infty}^2 + \frac{C\nu}{h^2}\right)\nu\|\nabla\ve{u}_{n+1}\|^2 \, ,
\end{align*}
which can be rewritten as
\begin{align*}
\frac{8}{7}\tau\beta_{n+1}\psi_{n+1}^2\|\nabla p_{n+1}\|^2 &\leq  4\tau\nu\|\nabla\ve{u}_{n+1}\|^2\\
&= \Phi_n-\Phi_{n+1} -  \bigg(\frac{\delta^2\psi_{n+1}}{\sqrt{\alpha}}\bigg)^2 -  \|\delta^2\ve{u}_{n+1}\|^2 - 4\tau\nu\gamma\|\nabla\cdot\ve{u}_{n+1}\|^2  ,
\end{align*}
with $\beta_{n+1}$ defined in \eqref{beta}. Finally, adding from $n=1$ to $n=N-1$ completes the proof.
}

\begin{remark}
    Very similar stability estimates can be shown for the first-order version of the method, which can be used to initialise the second-order scheme. 
\end{remark}

\section{Comparison with some classical methods}\label{sec_comparison}
Since the framework proposed here works as both a fractional-step scheme and a pressure stabilisation method, it is useful to compare the stability results above with those of classical methods that either decouple or stabilise the pressure. 

\subsection{Pressure-stabilised Petrov--Galerkin (PSPG)}
In the pressure estimate \eqref{stabilityP}, the factor $\beta_{n+1}$ resembles the parameter used in PSPG and similar stabilisation methods \cite{Codina1997} (although here it appears only in the analysis, not in the construction of the method). In other words, the present scheme works as a pressure stabilisation technique, circumventing the discrete inf-sup requirement. Differently from PSPG, however, we have here the velocity estimate \eqref{stabilityU} \textsl{decoupled} from the pressure stability, which may be seen as a form of pressure robustness (see Ref.~\cite{John2016}, section 4.6.2). Another key difference is that, with PSPG and similar methods, treating the pressure explicitly induces a parabolic CFL condition \cite{deFrutos2018}, so they should generally \textsl{not} be used as fractional-step schemes. Furthermore, the present method does not contain the term $\partial_t\ve{u}$ in the pressure equation; this is an advantage because the need to control the discrete acceleration in the PSPG residual can cause instabilities for small time-step sizes \cite{Burman2011}.

\subsection{Incremental pressure correction}
Incremental pressure-correction methods are the most popular fractional-step schemes. The \textsl{rotational} variant, which is the most accurate one, requires inf-sup-compatible pairs to guarantee pressure stability \cite{Guermond2006}. The classical (non-rotational) version, however, does have some added pressure stability due to the temporal discretisation, since it provides control over $\tau\sum_{n=1}^{N}\tau\|\nabla p_{n}-\nabla p_{n-1}\|^2$. However, as $\tau\rightarrow 0$ this term vanishes, which creates numerical instabilities \cite{Codina2001}. The pressure stability of the method proposed here, on the other hand, is proved in a much stronger norm that does not degrate for small $\tau$.

\section{Numerical examples}\label{sec_examples}
This section will test the accuracy and the stability of the new consistent scheme. All the examples use triangular finite elements with first-order Lagrangian interpolation for the pressure, paired with either quadratic or linear velocity spaces. The tests are performed at low to moderate Reynolds numbers $\mathrm{Re}\in [10,400]$, since convective stabilisation is not a topic under investigation here. The examples have no body forces ($\ve{f}=\ve{0}$).
	
\subsection{Temporal convergence test}
\begin{figure}[h]
    \centering
		\includegraphics[trim = 30 10 60 5,clip, width = .88\textwidth]{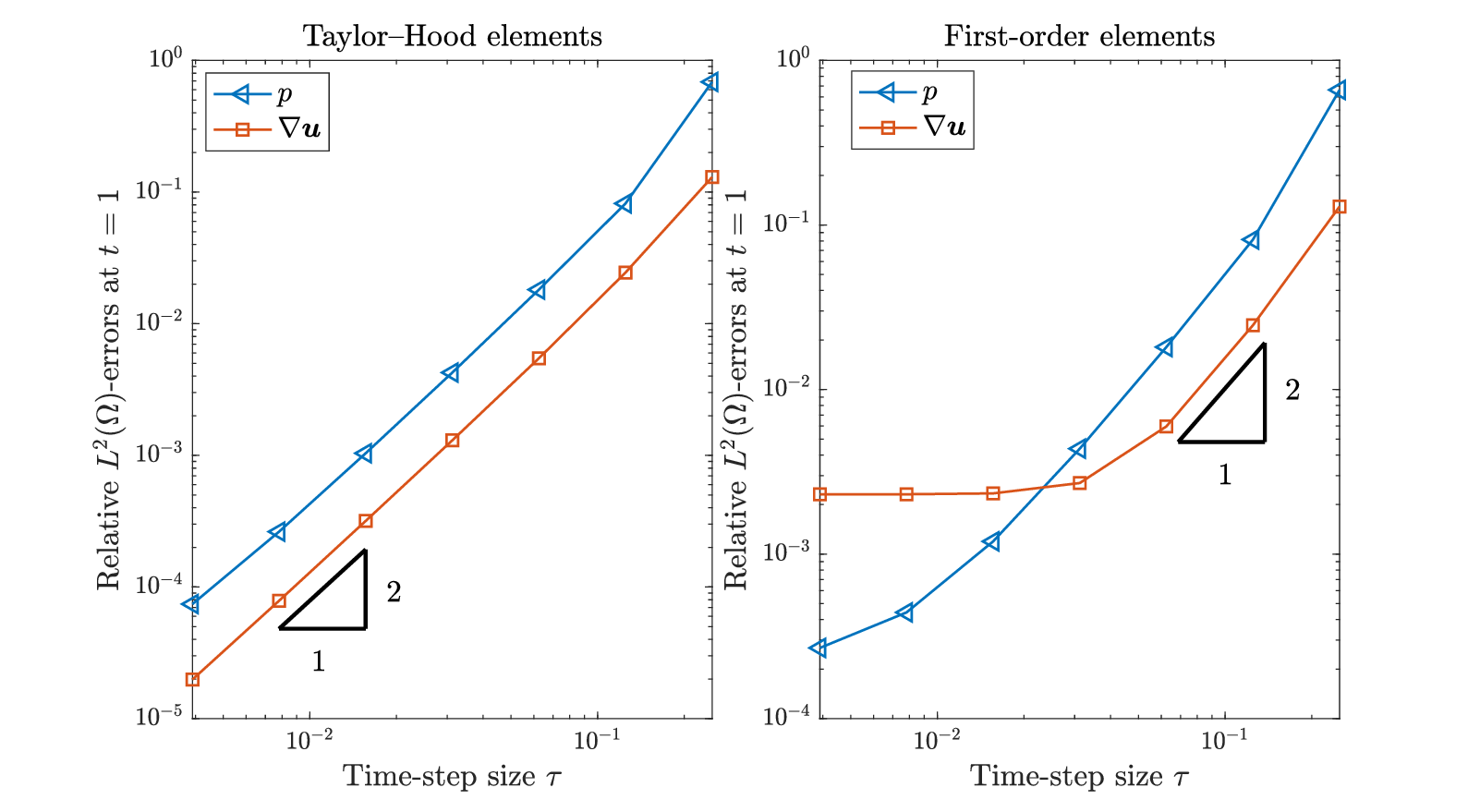}
    \caption{Taylor--Green vortex: temporal convergence considering a fine mesh with either Taylor--Hood (left) or equal-order (right) elements.}
    \label{TG}
\end{figure}

The first numerical example considers the Taylor-Green flow in $\Omega=(0,1)^2$, which has the analytical solution
	\begin{align*}
		\ve{u}  &= \mathrm{e}^{-2\pi^2\nu t}
		\begin{pmatrix}
			\hphantom{-}\sin\pi x\cos\pi y\\
			-\sin\pi y\cos\pi x
		\end{pmatrix}, \ \ p= \frac{\cos^2\pi x - \sin^2\pi y}{2}\mathrm{e}^{-4\pi^2\nu t} \, ,
	\end{align*}
describing a decaying vortex. This test case uses $\nu=0.1$, $\gamma=10$ and $\alpha=1$, and the Dirichlet boundary conditions for $\ve{u}$ are computed from the exact analytical expression. To verify that the discretisation method reaches full second-order convergence in time, a very fine mesh with 640,000 uniform elements is used. The time-step size starts at $\tau=0.25$ and is halved at each refinement level. The final time is $T=1$, at which point the $L^2(\Omega)$ errors for $p_{n+1}$ and $\nabla\ve{u}_{n+1}$ are computed. The results in Figure \ref{TG} confirm second-order temporal convergence for both Taylor--Hood and equal-order elements, with the latter breaking down earlier due to the larger spatial error. 

\subsection{Stability test}
\begin{figure}[h]
    \centering
		\includegraphics[trim = 10 10 10 10,clip, width = .7\textwidth]{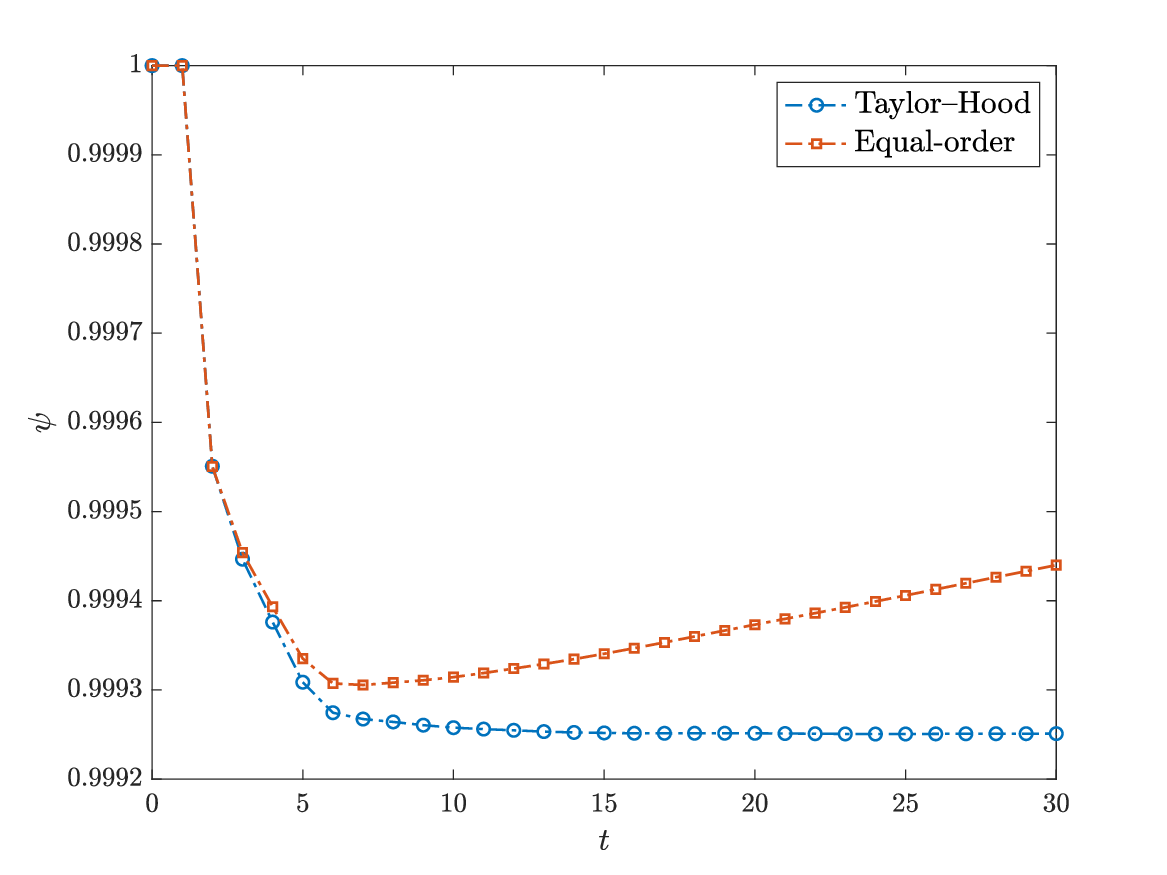}
    \caption{Lid-driven cavity flow: evolution of the SAV considering different types of elements.}
    \label{cavitySAV}
\end{figure}

The next example aims to put the splitting scheme to an extreme stability test by using a fine mesh and a very large $\tau$. The benchmark chosen is the lid-driven cavity flow at $\mathrm{Re}=400$. In $\Omega=(0,1)^2$, the fluid is initially at rest ($\ve{u}_0=\ve{0}$), the lateral and bottom walls are no-slip boundaries, and the lid velocity is $\ve{u}|_{y=1} = (u(x,t),0)^{\top}$, where
\begin{align*}
 &u(x,t) = (1-\mathrm{e}^{-3t})\left[1-\left(\frac{1-\cos\theta(x)}{2}\right)^2\right],
 \\
 &\theta(x) = \frac{\pi}{4l}\big(\big|2l + |2x-1| -1\big| + |2x-1| +2l-1\big)\, , \ \ l=0.1\, .
\end{align*}
The temporal part of the expression ramps up from 0 to 1, while the spatial part creates a smooth transition zone of length $l$ on each side of the lid, to regularise the corner discontinuity \cite{deFrutos2016}. The fluid and discretisation parameters are $\nu=0.0025$, $\gamma=100$ and $\alpha=0.1$. The mesh is as in the previous example, with either equal-order or Taylor--Hood elements. The time-step size is $\tau=1$ and the simulations run up to $t=30$, which is sufficiently long to reach a steady-state flow. The present results are compared to a reference stationary solution \cite{Ghia1982} with respect to the location of the primary vortex; the comparison, shown in Table \ref{table_lidVortex}, reveals very good agreement. Figure \ref{cavitySAV} depicts the evolution of the auxiliary variable $\psi$, which stays close to 1 throughout the simulation, as desired. Notice that the ratio $\tau\|\ve{u}\|_{\infty}/h$ for this example reaches $\mathcal{O}(10^3)$, which would be prohibitive under even a standard, hyperbolic CFL condition. 

\begin{table}[h]
 \centering
 \caption{Lid-driven cavity flow: steady-state location $(x,y)$ of the primary vortex.} 
 {\begin{tabular}{|c|c|c|}
    \hline
   \textbf{\citet{Ghia1982}}   & \textbf{Present (Taylor--Hood)} & \textbf{Present (equal-order)}\\
   \hline
   $(0.5547,0.6055)$   & $(0.555,0.605)$ & $(0.555,0.605)$ \\
    \hline
 \end{tabular}}
  \label{table_lidVortex} 
\end{table}

\subsection{Vortex shedding benchmark}    
\begin{figure}[h!]
    \centering
		\includegraphics[trim = 50 65 50 47,clip, width = .75\textwidth]{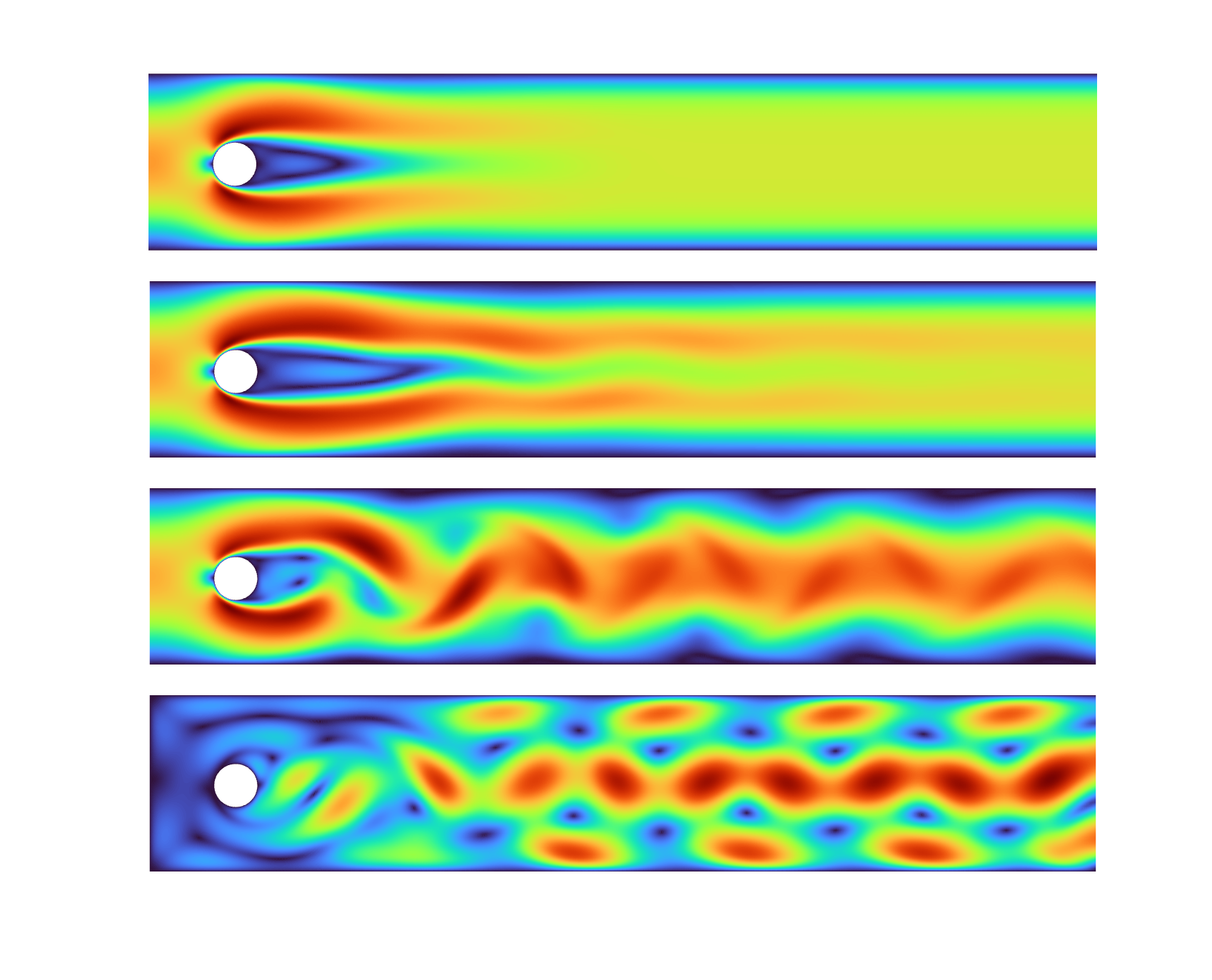}
    \caption{Cylinder flow benchmark: velocity magnitude fields for $t=2,4,6,8$.}
    \label{TurekShedding}
\end{figure}

The final numerical example aims to assess the new splitting method in the presence of high-frequency dynamics and a more general boundary condition setting. The example considered is the cylinder flow benchmark by \citet{Turek1996}. The setup, which has a time-periodic parabolic inflow and an open outflow, is classical and thoroughly described in the reference (see the geometry in Figure \ref{TurekShedding}). The simulation parameters are $\alpha=0.1$, $\gamma=1000$ and $\tau=0.0025$, using a structured mesh with 128,000 Taylor--Hood elements. The Reynolds number at peak inflow is $\mathrm{Re} = (2R) u_{\mathrm{max}}/\nu = 0.1\times 1/0.001 = 100$. 

\begin{figure}[h!]
    \centering
		\includegraphics[trim = 35 48 60 50,clip, width = .82\textwidth]{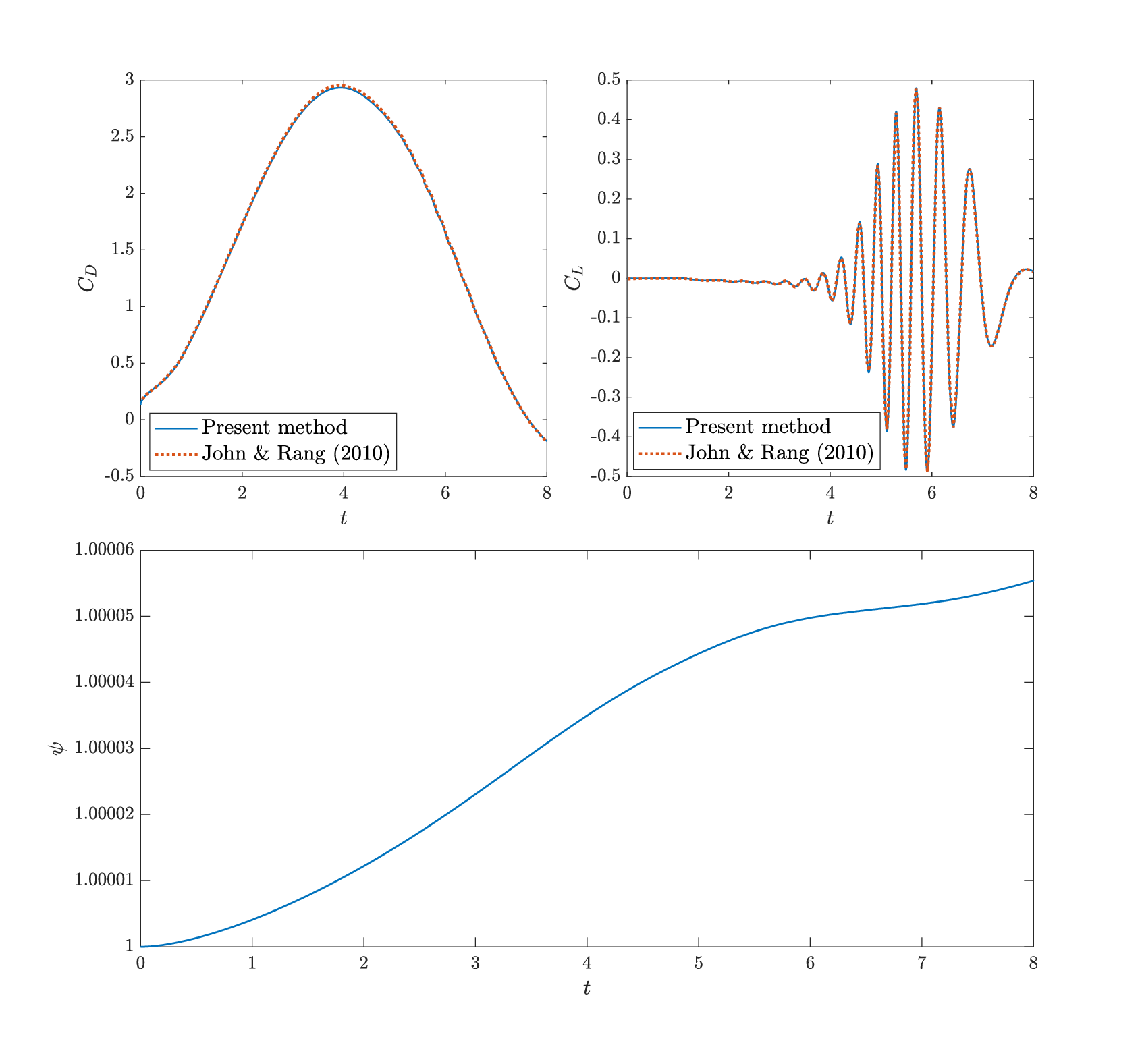}
    \caption{Cylinder flow benchmark: temporal evolution of the drag coefficient $C_D$, the lift coefficient $C_L$, and the auxiliary variable $\psi$.}
    \label{TurekCdCl}
\end{figure}

The velocity plots in Figure \ref{TurekShedding} show the evolution from a quasi-symmetric flow field to a vortex shedding scenario. The main benchmark quantities are the cylinder's lift and drag coefficients, $C_L$ and $C_D$, for which Figure \ref{TurekCdCl} shows excellent agreement between the present solution and reference results \cite{John2010}. The peak coefficients obtained here were ${C_D}_{\mathrm{max}}\approx 2.9327$ and ${C_L}_{\mathrm{max}}\approx 0.47951$, against ${C_D}_{\mathrm{max}}\approx 2.9509$ and ${C_L}_{\mathrm{max}}\approx 0.47787$ reported by \citet{John2010}. As in the previous example, the SAV $\psi$ remains very close to unity throughout the simulation.

\newpage	
\section{Concluding remarks}\label{sec_Conclusion}
	This article has presented what is, to the best of the author's knowledge, the first provenly stable consistent splitting scheme using standard finite elements. To guarantee stability under standard assumptions, an auxiliary time-dependent unknown is introduced, in the spirit of recent SAV methods. The absence of splitting errors and spurious boundary conditions allows the method to achieve the full order of the temporal discretisation, which is usually not possible in standard fractional-step methods---especially not in the presence of open (Neumann) boundaries. The resulting scheme can be implemented very efficiently, requiring only the solution of linear advection-diffusion-reaction equations and a Poisson problem. Unconditional temporal stability has been proved for a prototypical second-order scheme under standard assumptions on the finite element discretisation. The numerical examples confirm second-order convergence and unconditional stability through various problems of increasing complexity, which reach extreme CFL numbers surpassing 1000. Some questions and challenges remain open, such as a fully discrete error analysis for the first- and second-order BDF methods. Especially relevant would be analysing stability for a BDF3 scheme, whose order cannot be reached by pressure-correction methods. 
    
	\section*{Acknowledgments}
	The author acknowledges funding by the Federal Ministry of Education and Research (BMBF) and the Ministry of Culture and Science of the German State of North Rhine-Westphalia (MKW) under the Excellence Strategy of the Federal Government and the Länder.

   \section*{Data availability}
   No data was used for or generated from the research described in this article.

\section*{Declarations}
\textbf{Conflict of interest}. The author declares that there are no conflicts of interest.
	
	\bibliographystyle{unsrtnat} 
	\bibliography{Preprint}%
	
\end{document}